\documentclass{amsart}
\usepackage[hidelinks]{hyperref}

\usepackage{amssymb}

\def\dom{\mathop{\mathrm{Dom}}\nolimits}

\def\str#1{\mathbf {#1}}
\def\Fraisse{Fra\"{\i}ss\' e}

\begin{document}


\title{Extending partial automorphisms of $n$-partite tournaments}


\author[J. Hubi\v cka]{Jan Hubi\v cka}
\address{Department of Applied Mathematics (KAM), Charles University, Prague, Czech Republic}
\email{hubicka@kam.mff.cuni.cz}


\author[C. Jahel]{Colin Jahel}
\address{Universit\'e Paris Diderot, Institut de Math\'ematiques de Jussieu-Paris Rive Gauche}
\email{colin.jahel@imj-prg.fr}

\author[M. Kone\v cn\'y]{Mat\v ej Kone\v cn\'y}
\address{Department of Applied Mathematics (KAM), Charles University, Prague, Czech Republic}
\email{matej@kam.mff.cuni.cz}


\author[M. Sabok]{Marcin Sabok}
\address{Department of Mathematics and Statistics, McGill University, 805, Sherbrooke Street West Montreal, Quebec, Canada H3A 2K6, and
Institute of Mathematics, Polish Academy of Sciences, Sniadeckich 8, 00-655 Warszawa, Poland}
\email{marcin.sabok@mcgill.ca}


\thanks{J. Hubi\v cka and M. Kone\v cn\'y are supported by project 18-13685Y of the Czech Science Foundation (GA\v CR) and by Charles University projects Progres Q48 and GA UK No 378119 respectively. M. Sabok is supported by project FRQNT (Fonds de recherche du Quebec) grant Nouveaux chercheurs 2018-NC-205427. C. Jahel is partially supported by the ANR project AGRUME
(ANR-17-CE40-0026).}





\begin{abstract}
We prove that for every $n\geq 2$ the class of all finite $n$-partite tournaments (orientations of complete $n$-partite graphs) has the extension property for partial automorphisms, that is, for every finite $n$-partite tournament $G$ there is a finite $n$-partite tournament $H$ such that every isomorphism of induced subgraphs of $G$ extends to an automorphism of $H$. Our constructions are purely combinatorial (whereas many earlier EPPA results use deep results from group theory) and extend to other classes such as the class of all finite semi-generic tournaments.
\end{abstract}

\maketitle


\theoremstyle{plain}
\newtheorem{theorem}{Theorem}
\newtheorem{question}{Question}

\theoremstyle{remark}
\newtheorem*{remark}{Remark}

We say that a class of finite structures $\mathcal C$ has the \emph{extension property for partial automorphisms} (\emph{EPPA}, also called the \emph{Hrushovski property}) if for every $\str A\in \mathcal C$ there is an \emph{EPPA-witness} $\str B\in \mathcal C$, that is, a structure containing $\str A$ as an (induced) substructure such that every isomorphism of substructures of $\str A$ (a \emph{partial automorphism} of $\str A$) extends to an automorphism of $\str B$.

In 1992 Hrushovski~\cite{hrushovski1992} established that the class of all finite graphs has EPPA. This results, besides being a beautiful combinatorial property by itself, was used by Hodges, Hodkinson, Lascar and Shelah to study the automorphism group of the countable random graph~\cite{hodges1993b} and found many additional applications~\cite{Siniora2}. After this, the quest of identifying new classes of structures with EPPA continued. 
In 2000, Herwig and Lascar~\cite{herwig2000} combined advanced techniques from combinatorics, model theory and group theory and gave a structural condition which was later used to identify many additional examples of classes having EPPA (such as various classes of metric spaces~\cite{solecki2005,Hubicka2018metricEPPA,Aranda2017,Konecny2019a}).

Herwig and Lascar also asked whether
the class of all finite tournaments has EPPA.  This is probably the simplest class that can not be handled by the Herwig--Lascar theorem and is of special interest for group-theorists. While this question remains among the most important open questions in the area, a negative result was recently obtained by Huang, Pawliuk, Sabok and Wise that the class of certain generalised tournaments does not have EPPA~\cite{Sabok}. On the other hand, EPPA was shown for the class of two-graphs by Evans, Hubi\v cka, Kone\v cn\'y and Ne\v{s}et\v{r}il~\cite{eppatwographs}.  This class behaves similarly to tournaments from the perspective of the Herwig--Lascar theorem. We further develop the techniques used in~\cite{eppatwographs} to show:
\begin{theorem}
\label{thm:partite}
For every $n\geq 2$, the class of all finite $n$-partite tournaments has the extension property for partial automorphisms.
\end{theorem}
Here, a directed graph $G=(V,E)$ is an \emph{$n$-partite tournament} if its vertex set can
be partitioned into sets $V_1\cup V_2 \cup\ldots \cup V_n=V$ (called \emph{parts}) such that
every pair of vertices in different parts is connected by exactly one directed edge and
there are no edges between vertices in the same part.

\medskip

Unlike many earlier proofs of various EPPA results, the proof of Theorem~\ref{thm:partite} is purely combinatorial and fully self-contained. Developing the idea of valuation functions introduced by Hodkinson and Otto~\cite{hodkinson2003} (which has also been used to strengthen the Hodkinson--Otto result to classes with unary functions~\cite{Evans3} and applied to give a combinatorial proof of EPPA for metric spaces~\cite{Hubicka2018metricEPPA} and a strengthening of the Herwig--Lascar theorem~\cite{Hubicka2018EPPA}), we construct highly symmetric $n$-partite tournaments to serve as EPPA-witnesses.

EPPA is an important concept in model theory, because it can be seen as a strengthening of the \emph{amalgamation property} (see e.g.~\cite{Macpherson2011}).
By the \Fraisse{} theorem, classes with this property in turn correspond to the (most often countably infinite) \emph{homogeneous} structures, that is, structures such that every isomorphism of their finite substructures extends to an automorphism (i.e., they are infinite EPPA-witnesses for themselves).

The goal of the Lachlan--Cherlin classification programme is to classify homogeneous structures (subject to extra conditions).
In particular, Cherlin's classification of homogeneous directed graphs~\cite{Cherlin1998} gives a complete list of amalgamation classes
of directed graphs. This relatively elaborate list was analysed by Pawliuk and Soki\'c~\cite{PawliukSokic16} and up three exceptions it was possible to decide for each class whether it has EPPA or not.
The remaining three cases are the classes of all $n$-partite tournaments for $n\geq 2$, the class of all semi-generic tournaments and the classes of
all oriented graphs omitting an independent set of size $k$ for $k\geq 2$ (for $k=2$ this is the class of all tournaments).

A more involved adaptation of the techniques presented here also gives:
\begin{theorem}
\label{thm:semigeneric}
The class of all semi-generic tournaments has EPPA.
\end{theorem}
(An $n$-partite tournament $G$ with parts $V_1,V_2,\cdots, V_n$ is \emph{semi-generic} if for every $i\neq j$ and every choice of $a\neq b\in V_i$ and $c\neq d\in V_j$
it holds that the number of edges directed from $a$ or $b$ to $c$ or $d$ is even.) This will appear in~\cite{HubickaSemigeneric}.

\section*{Proof of Theorem~\ref{thm:partite}}

Fix $n\geq 2$ and a finite $n$-partite tournament $G=(V,E)$ with parts
$V_1,V_2,\ldots,V_n$. Without loss of generality we can assume the following:
\begin{enumerate}
\item $V=\{1,2,\ldots, k\}$,
\item for every $x\in V_i$ and every $y\in V_j$ it holds that $x<y$ whenever $i<j$, and 
\item $|V_1|=|V_2|=\ldots=|V_n|$.
\end{enumerate}

We give an explicit construction of an $n$-partite tournament $H$ which is an EPPA-witness for $G$.
Given $x\in V$, a \emph{valuation function for $x$}
is a function $\chi_x\colon V\setminus\{x\}\to \{0,1\}$ with the property that $\chi_x(y)=0$ for all $y$ which are in the same part as $x$.
For a valuation function $\chi_x$ for $x$ we define its \emph{projection} as $\pi(\chi_x)=x$.
For clarity, whenever we use $\chi_x$ for a valuation function, it means that $\pi(\chi_x)=x$.

The directed graph $H=(V',E')$ is constructed as follows:
\begin{enumerate}
 \item The vertex set $V'$ of $H$ consists of all valuation functions for all vertices $x\in V$, and
 \item $\chi_x,\chi_y\in V'$ are adjacent if and only if $x$ and $y$ belong to different parts of $G$. The edge
is oriented from $\chi_x$ to $\chi_y$ if and only if one of the following is satisfied:
\begin{itemize}
\item $x > y, \text{ and }\chi_x(y)\neq \chi_y(x)$, or
\item $x < y, \text{ and }\chi_x(y) = \chi_y(x)$.
\end{itemize}
Otherwise the edge is oriented from $\chi_y$ to $\chi_x$.
\end{enumerate}
Clearly, $H$ is an $n$-partite tournament with parts $V'_i=\{\chi_x\in V':x\in V_i\}$.

\medskip

Next we construct an embedding $\psi\colon G\to H$. Given a vertex $x\in V$ we put $\psi(x) = \chi_x$, where $\chi_x$ is a function $V\setminus\{x\}\to \mathbb \{0,1\}$ such that
$$\chi_x(y)=
\begin{cases}
1 & \text{ if } y<x \text{ and there is an edge directed from $x$ to $y$ in $G$}\\
0 & \text{ otherwise}.
\end{cases}$$
It is easy to verify that $\psi$ is indeed an embedding of $G$ into $H$. Put $G'=\psi(G)$.

\medskip

We now show that $H$ extends all partial automorphisms of $G'$.
Fix a partial automorphism $\varphi\colon G'\to G'$.  By the projection,
$\varphi$ induces a partial permutation of $V$, which in turn induces a partial permutation 
$\iota$ of parts $V_1, V_2,\ldots, V_n$. Let $\hat\iota$ be an extension of $\iota$ to a full permutation of parts.
Denote by $\hat\varphi$ a permutation of $V$ which extends the permutation induced on $V$
by $\varphi$ and which furthermore agrees with $\hat\iota$ (such a permutation exists
as we assume all parts to have the same size).

For every pair of vertices $x<y$ of $G$ we define a function $F_{\{x,y\}}\colon \{x,y\}\to \{0,1\}$ as follows:
\begin{enumerate}
 \item $F_{\{x,y\}}(x)=F_{\{x,y\}}(y)=0$ if $x$ and $y$ belong to the same part.
 \item If $x$ and $y$ belong to different parts and $\hat\varphi(x)<\hat\varphi(y)$ we put
 \begin{enumerate}
   \item $F_{\{x,y\}}(x)=F_{\{x,y\}}(y)=1$ if
\begin{enumerate}
\item $\psi(x)\in \dom(\varphi)$ and $\psi(x)(y)\neq \psi(\hat\varphi(x))(\hat\varphi(y))$ or
\item  $\psi(y)\in \dom(\varphi)$ and $\psi(y)(x)\neq \psi(\hat\varphi(y))(\hat\varphi(x))$,
\end{enumerate}
   \item $F_{\{x,y\}}(x)=F_{\{x,y\}}(y)=0$ otherwise.
 \end{enumerate}
 \item if $x$ and $y$ belong to different parts and $\hat\varphi(x)>\hat\varphi(y)$ we put
 \begin{enumerate}
   \item $F_{\{x,y\}}(x)=1$ and $F_{\{x,y\}}(y)=0$ if
\begin{enumerate}\item $\psi(x)\in \dom(\varphi)$ and $\psi(x)(y)\neq \psi(\hat\varphi(x))(\hat\varphi(y))$ or
\item $\psi(y)\in \dom(\varphi)$ and $\psi(y)(x)= \psi(\hat\varphi(y))(\hat\varphi(x))$,
\end{enumerate}
   \item $F_{\{x,y\}}(x)=0$ and $F_{\{x,y\}}(y)=1$ otherwise.
 \end{enumerate}
\end{enumerate}

Now we are ready to describe a function $\theta\colon V'\to V'$ which we will verify to be an automorphism of $H'$ extending $\varphi$. We put $\theta(\chi_x)=\chi'$, where $\chi'$ is a function $V\setminus \{\hat\varphi(x)\}\to \{0,1\}$ defined as:
$$
\chi'(\hat\varphi(y))=
\begin{cases}
  \chi_x(y) &  \text{if }y\neq x\text{ and } F_{\{x,y\}}(x)=0 \\
  1-\chi_x(y) &  \text{if }y\neq x\text{ and } F_{\{x,y\}}(x)=1\\
  \text{undefined} & \text{if }y=x.
\end{cases}
$$

\medskip

It is easy to verify that $\theta$ is one-to-one, because one can construct its inverse. To verify that $\theta$ extends $\varphi$, we want to show that for every $\chi_x\in \dom(\varphi)$ it holds that
$\varphi(\chi_x)(y)=\psi(\hat\varphi(x))(y)$. This can be checked directly using the definition of functions $F_{\{x,y\}}$.

Finally we need to prove that $\theta$ is an automorphism of $H$. Clearly $\theta$ preserves non-edges. And the functions $F_{\{x,y\}}$ were chosen precisely so that $\theta$ also preserves the directions of edges.

\begin{remark}
Note that the constructed $H$ is an EPPA-witness for all $n$-partite tournaments on at most $|G|$ vertices.
\end{remark}

\section*{Conclusion and remarks}
Our proof uses the method of valuation functions which has multiple additional applications mentioned above.

In general, to use the method, one needs to first fix $\str A$ whose partial automorphism we will extend and to understand what these partial automorphisms look like. Then one defines a suitable notion of valuation functions such that ``each vertex can change its neighbourhood arbitrarily''. A correct notion of flips then follows from this interpretation and it turns out that putting $\str B$ to simply be all valuation functions for all vertices of $\str A$ often works. See~\cite{Hubicka2018EPPA,Hubicka2018metricEPPA,eppatwographs} for examples.

The proof of Theorem~\ref{thm:semigeneric} also goes in this direction, however, one has to give valuation functions not only to vertices, but also to parts.

\begin{remark}
As noted above, every class with EPPA is an amalgamation class and the study has thus so far been limited to them. However, it might be reasonable to refine the problem and ask ``Given a class $\mathcal C$, for which structures $\str A\in \mathcal C$ is there an EPPA-witness $\str B\in\mathcal C$?'' It remains to be seen if this leads to some interesting examples.
\end{remark}

\begin{question}
The EPPA witnesses constructed in this paper have $\mathcal O(n2^n)$ vertices, the valuation functions construction for graphs also gives this number. There are, however, no non-trivial lower bounds for the size of EPPA-witnesses for structures on at most $n$ vertices. Is it possible to close this gap?
\end{question}


\section*{Acknowledgement}
This research was started during the workshop \textit{Unifying Themes in Ramsey Theory} which took place in Banff in November 2018.


\bibliographystyle{abbrv}
\bibliography{ramsey.bib}





\end{document}